\documentclass{article}

\def \R {{\mathbf {R}}}
\def \B  {\cal B}

\def\uu{\bigsqcup}
\def\eps{\varepsilon}

\usepackage[T1]{fontenc} % РёРЅРѕРіРґР° T1 РІРјРµСЃС‚Рѕ T2A
\usepackage[cp1251]{inputenc} 
\usepackage[russian]{babel}

\begin{document}

  \rm \Large

\begin{center}   
{ \LARGE \bf
 Typical infinite mixing automorphisms are  rank-one  
\vspace{5mm}

Valery V. Ryzhikov
}
\end{center} 

\vspace{5mm} \rm
Bashtanov proved that generic mixing automorphisms of probability space with respect to the Alpern-Tikhonov metric had rank one. Using Sidon's constructions, we show that generic infinite mixing automorphisms also are  rank-one.

\section{Introduction}
Alpern and Tikhonov \cite{A},\cite{T} introduced a complete metric on the space of mixing actions. Bashtanov \cite{Ba} proved that the class of rank-one actions is generic in the space ${\bf Mix}$ of invertible mixing transformations with invariant probability measure. We will establish the same  for the space ${\bf Mix_\infty}=\{T:\, T^n\to_w 0\}$ of transformations of $T$ with sigma-finite invariant measure. We  prove that the conjugacy class of the Sidon transformation of rank one (see \cite{24},
\cite{24a}) is dense in ${\bf Mix_\infty}$. Since the class of rank-one transformations is a $G_\delta$-set, taking into account the density of the latter, we obtain the genericity of this class. A consequence of this  and  the result of \cite{15} is the typicalness of the minimality of the centralizer (the tipical transformation $T$ commutes only with its powers).

\vspace{2mm}
The group ${\bf Aut_\infty }$ of all atomorphisms of the space
$(X,\B,\mu)$ isomorphic to the line $\R$ with the standard Lebesgue measure
is equipped with the Halmos metric $\rho$:
$$ \rho(S,T)=\sum_i 2^{-i}\left(\mu(SA_i\Delta TA_i)+\mu(S^{-1}A_i\Delta T^{-1}A_i)\right),$$
where $\{A_i\}$ is some fixed family dense in the family of all measurable sets of finite measure.
The space $({\bf Aut_\infty},\rho)$ is a complete separable metric space.
It is easy to show that in $({\bf Aut_\infty},\rho)$ the so-called mixing transformations form a set of
the first category (countable intersection of nowhere dense sets).
If we strengthen the Halmos metric to the Alpern-Tikhonov metric, we obtain a complete metric space ${\bf Mix}_\infty$.
Which invariants (i.e. properties preserved under conjugation in $\bf Aut_\infty$) are typical? A property is called typical if
a dense $G_\delta$-set of transformations has this property.

\vspace{3mm}
\bf Spaces $\bf {\bf Mix}$, ${\bf Mix}_\infty$. \rm
Following \cite{T}, we define the metric $d_w$ on ${\bf Aut}(\mu)$:
$$ d_w(S,T)=\sum_{i,j=1}^\infty 2^{-i-j}\left|\mu(SA_i\cap A_j)-\mu(TA_i\cap A_j)\right|.$$
On the set of all mixing automorphisms {\bf Mix} is defined the metric $r$:
$$ r(S,T)= \rho(S,T) + \sup_{n>0}d_w(S^n,T^n).$$
In \cite{T} it is proved that the metric space $({\bf Mix},r)$ is complete and separable.

Consider the sigma-finite measure $\mu_\infty$ as an invariant and, without changing the formulas, introduce the same metric $r$ on the space of mixing infinite automorphisms:
$${\bf Mix}_\infty =\{ T\in {\bf Aut}(\mu_\infty)\, : \, T^n\to_w 0 \},$$
where $\{A_i\}$ is a fixed sequence of sets,
dense in the family of all sets of finite measure.

Fundamental in $({\bf Mix_\infty}, r)$  a sequence $T_k$ converges to some $T$ by virtue of the completeness of the Halmos metric
in the space ${\bf Aut}(\mu_\infty)$. It is easy to show that this
$T_k\to T$ converges to $T$  in the metric $r$ as well.

From the known methods, it can be shown that 

\it for generic  transformations in ${\bf Mix}_\infty$ that  spectra of their  symmetric powers are simple. \rm

\vspace{3mm}
It  follows from the results of \cite{24a} that

\vspace{3mm}
\it a typical transformation in ${\bf Mix}_\infty$
is not isomorphic to its inverse. \rm

\section{Density of the conjugacy classes of the  Sidon automorphisms}
Let us recall  necessary definitions.

\vspace{2mm}
\bf Construction of rank one. \rm
Fix a natural number $h_1$, a sequence $r_j\to\infty$ (parameter $r_j$ is the number of columns into which the tower of stage $j$ is cut) and a sequence of integer vectors (parameters of spacers)
$$ \bar s_j=(s_j(1), s_j(2),\dots, s_j(r_j-1),s_j(r_j)), \ s_(i)\geq 0.$$
At step $j=1$ a set of disjoint half-intervals
$E_1, SE_1, S^2E_1,\dots, S^{h_1-1}E_1$ of the same measure is given. 

Let at step $j$ a system of non-intersecting half-intervals
$$E_j, SE_j, S^2E_j,\dots, S^{h_j-1}E_j,$$ be defined, and on the half-intervals $E_j, SE_j,\dots, S^{h_j-2}E_j$
the transformation $S$ is the parallel translation. Such a set of half-intervals is called the tower of step $j$, their union is denoted by $X_j$ and is also called a tower.

We represent $E_j$ as a disjoint union of $r_j$ half-intervals
$$E_j^1,E_j^2E_j^3,\dots E_j^{r_j}$$ of the same measure (length).
For each $i=1,2,\dots, r_j$, we consider the so-called column
$$E_j^i, SE_j^i ,S^2 E_j^i,\dots, S^{h_j-1}E_j^i.$$
We denote the union of these half-intervals by $X_{i,j}$, which we will also call a column.

To each column $X_{i,j}$ we add $s_j(i)$ disjoint half-intervals of the same measure as $E_j^i$, obtaining a set
$$E_j^i, SE_j^i, S^2 E_j^i,\dots, S^{h_j-1}E_j^i, S^{h_j}E_j^i, S^{h_j+1}E_j^i, \dots, S^{h_j+s_j(i)-1}E_j^i$$
(all these sets are disjoint).
Denoting $E_{j+1}= E^1_j$, for $i<r_j$ we set
$$S^{h_j+s_j(i)}E_j^i = E_j^{i+1}.$$
The set of superstructured columns is from now on considered as a tower of stage $j+1$, consisting of half-intervals
$$E_{j+1}, SE_{j+1}, S^2 E_{j+1},\dots, S^{h_{j+1}-1}E_{j+1},$$
where
$$ h_{j+1} =h_jr_j +\sum_{i=1}^{r_j}s_j(i).$$

The partial definition of the transformation $S$ at stage $j$ is preserved at all subsequent stages. As a result, on the space $X=\cup_j X_j$ there is defined an invertible transformation $S:X\to X$ that preserves the standard Lebesgue measure on $X$.

\vspace{2mm}
\bf Sidon constructions. \rm Let a rank one construction $S$  have the following property: \it the intersection
$X_j\cap S^mX_j$ for $h_{j}<m\leq h_{j+1}$ can be contained
only in one of the columns $X_{i,j}$ of the tower $X_j$. \rm Such a construction  is called  Sidon. \rm The measure of  its space $X$  is infinite. If $r_j\to\infty$, then Sidon $S$ is mixing.

\vspace{2mm}
\bf Theorem 1. \it Let $T$ be a mixing transformation, $S$ a Sidon mixing transformation of rank one. For every $\eps>0$ there exists a conjugation $R$ such that $ r(R^{-1}SR,T)<\eps$.\rm

\vspace{2mm}
\bf Lemma 2. \it Let $T$ be a mixing transformation, $A$ a set of finite measure. The phase space for every natural $H$ can be represented as $$X=\uu_{h>H}\uu_{i=0}^{h}T^iB_{h}, \ \ \ T^hB_{h}\cap A=\phi.$$ \rm

\vspace{2mm}
The proof of the lemma is simple, let it be  an exercise.

\vspace{2mm}
The proof of the theorem (in brief). 

Let us consider a part of the Kakutani tower:
$$U_N=\uu_{H<h<N}\uu_{i=0}^{h}T^iB_{h},$$
which almost absorbed the sets $A_1,\dots, A_k$, appearing in the definition of the metric $r$. For this, we choose a large $N$.

We take a very thin tower $X_j$ (so very large $j$), appearing in the construction of the Sidon transformation $S$, and wind on $U_N$ the tower $\tilde X_j$ corresponding to the conjugation
$ \tilde S=R^{-1}SR$ in  a way that the following conditions are satisfied:

(1) the conjugate transformation $\tilde S$ coincides with $T$ on $U_N$ minus the upper floors of $T^hB_{h}$ and minus a set of arbitrarily small measure;

(2) after the moving floor of the wound tower leaves the top of $U_N$ for the $k$-th time, it must travel outside $U_N$ for a time $s_k> 10 s_{k-1}$
($s_1> 10N$). This  Sidon winding
will ensure that the intersection measure is small:
$$\mu(\tilde S^m\tilde U_N\cap U_N)\ <\ N\mu(B_j), \ \ m<h_j.$$

For $m\geq h_j$, the inequality
$$\mu(\tilde S^m\tilde U_N\cap \tilde U_N)<\mu(\tilde U_N)/r_j$$ is satisfied by  the definition of the Sidon construction.

If $A=A_1\cup\dots \cup A_k\subset \tilde U_N$,
then
$$\sum_i^k 2^{-i}\left(\mu(SA_i\Delta TA_i)+\mu(S^{-1}A_i\Delta T^{-1}A_i)\right)=0$$
and
$$\sum_{i,j=1}^k 2^{-i-j}\left|\mu(SA_i\cap A_j)-\mu(TA_i\cap A_j)\right|=0.$$
In reality, for a large $N$ and a suitable Sidon winding, instead of equalities we get inequalities with a very small right-hand side.

So for $\eps>0$ there are  $j$, $N$ and $R$ such that at the intersection $\tilde X_j\cap U_N$ the transformations $\tilde S$ and $T$ almost coincide, hence,  $\tilde S$ and $T$ are close in the Halmos metric,
but their powers are uniformly close with respect to the $d_w$ metric.
Thus we obtain $r(R^{-1}SR,T)<\eps$ for a suitable $R$.

The set of rank-one transformations is a $G_\delta$ set
(use the proof from \cite{Ba}).
We have found that the class of mixing rank one transfomations  is dense
in $\bf Mix_\infty$.  Thus, it is typical.

\vspace{3mm}
In \cite{20} we used a similar trick, but the windings were random Ornstein 
and the spacers for flows and automorphisms of rank one were random Ornstein as well.

In the case of infinite measure, a pleasant opportunity arises to use large Sidon spacers, which gives  us  a  simpler and  constructive method. For other results on Sidon transformations, see in  \cite{24}.

\section{Типичность ранга один в  ${\bf Mix_\infty}$}
Альперн и Тихонов  \cite{A},\cite{T} ввели полную метрику на пространстве перемешивающих действий.   Баштанов \cite{Ba} доказал, что    класс действий ранга один  типичен в пространстве ${\bf Mix}$ обратимых перемешивающих преобразований с   инвариантной вероятностной мерой. Мы установим аналогичную теорему для   пространства ${\bf Mix_\infty}=\{T:\, T^n\to 0\}$  преобразований  $T$ с сигма-конечной инвариантной мерой.  Докажем, что  класс сопряженности  сидоновского преобразования ранга один (см. \cite{24}, 
\cite{24a})  плотен в  ${\bf Mix_\infty}$.  Так как класс преобразований ранга один является $G_\delta$-множеством, с учетом  плотности последнего получаем типичность этого класса. Следствием этой типичности с учетом  результата \cite{15}  является типичность минимальности централизатора (преобразование $T$ коммутирует только со своими степенями).

Группа  ${\bf Aut_\infty }$ всех атоморфизмов  пространства 
$(X,\B,\mu)$, изоморфного прямой $\R$ со стандартной  мерой  Лебега
оснащается метрикой Халмоша    $\rho$:
$$ \rho(S,T)=\sum_i 2^{-i}\left(\mu(SA_i\Delta TA_i)+\mu(S^{-1}A_i\Delta T^{-1}A_i)\right),$$
где $\{A_i\}$  -- некоторое  фиксированное семейство,  плотное в семействе всех  измеримых множеств конечной меры.
Пространство  $({\bf Aut_\infty},\rho)$  является полным сепарабельным пространством.
Несложно показать, что в  $({\bf Aut_\infty},\rho)$ так называемые перемешивающие преобразования образуют множство 
первой категории (счетное пересечение нигде не плотных множеств).
Если же усилить метрику Халмоша до метрики Альперна-Тихонова, получим полное метрическое пространство  ${\bf Mix}_\infty$). 
Возникает    вопрос:  какие инварианты, т.е. свойства, сохранящиеся при сопряжении в $\bf Aut_\infty$) являются типичными? Свойство называется типичным, если 
плотное $G_\delta$-множество преобразований   обладают этим свойством.

Например, с учетом  известных методов   можно показать, что 
в ${\bf Mix}_\infty$ типична однократность спектра всех симметрических степеней преобрзования.  Из результатов \cite{24a} легко вытекает, что 
\vspace{3mm}
\it типичное в ${\bf Mix}_\infty$ 
преобразование     не изоморфно своему   обратному. \rm

\vspace{3mm}
\bf Пространства $\bf {\bf  Mix}$, ${\bf  Mix}_\infty$. \rm
Cледуя \cite{T}, определим   на ${\bf Aut}(\mu)$ метрику  $d_w$:
$$ d_w(S,T)=\sum_{i,j} 2^{-i-j}\left|\mu(SA_i\cap A_j)-\mu(TA_i\cap A_j)\right|.$$
На множестве всех перемешиващих автоморфизмов {\bf  Mix} задаем  метрику $r$:
$$  r(S,T)= \rho(S,T) + \sup_{n>0}d_w(S^n,T^n).$$
В \cite{T} доказано, что метрическое пространство $({\bf  Mix},r)$  является полным и сепарабельным.

Рассмотрим в качестве инвариантной сигма-конечную меру $\mu_\infty$ и, не меняя формул, введем такую   же метрику $r$ на пространстве перемешивающих бесконечных автоморфизмов:
$${\bf  Mix}_\infty =\{ T\in {\bf Aut}(\mu_\infty), \, T^n\to_w 0 \},$$
где  $\{A_i\}$  --  фиксированная последовательность  множеств,  
плотная в семействе всех множеств конечной меры.

Фундаментальная в $({\bf  Mix_\infty}, r)$  последовательность $T_k$  сходится к некоторому $T$ в силу полноты метрики Халмоша
в пространстве ${\bf Aut}(\mu_\infty)$. Несложно показать, что в этом случае
$T_k\to T$  и в метрике $\rho$.

\bf Плотность класса сопряженности сидоновского автоморфизма. \rm

\vspace{2mm}
\bf   Конструкция ранга один. \rm

\vspace{2mm}
Фиксируем натуральное число $h_1$, последовательность  $r_j\to\infty$ ( параметр $r_j$ -- число колонн, на которые  разрезается башня этапа $j$) и  последовательность целочисленных векторов (параметров надстроек)   
$$ \bar s_j=(s_j(1), s_j(2),\dots, s_j(r_j-1),s_j(r_j)).$$  
На шаге $j=1$ по определению задан  набор полуинтервалов 
$E_1, SE_1, S^2E_1,\dots, S^{h_1-1}E_1$. Пусть на   шаге $j$  определена система   непересекающихся полуинтервалов 
$$E_j, SE_j, S^2E_j,\dots, S^{h_j-1}E_j,$$
причем на полуинтервалах $E_j, SE_j, \dots, S^{h_j-2}E_j$
пребразование $S$ является параллельным переносом. Такой набор   полуинтервалов  называется башней этапа $j$, их объединение обозначается через $X_j$ и тоже называется башней.

Представим   $E_j$ как дизъюнктное объединение  $r_j$ полуинтервалов 
$$E_j^1,E_j^2E_j^3,\dots E_j^{r_j}$$ одинаковой меры (длины).  
Для  каждого $i=1,2,\dots, r_j$ рассмотрим так называемую колонну  
$$E_j^i, SE_j^i ,S^2 E_j^i,\dots, S^{h_j-1}E_j^i.$$
Объдинение этих полуинтервалов  обозначим через $X_{i,j}$, его тоже будем называть колонной.

К каждой  колонне $X_{i,j}$ добавим  $s_j(i)$  непересекающихся полуинтервалов  той же меры, что у $E_j^i$, получая набор  
$$E_j^i, SE_j^i, S^2 E_j^i,\dots, S^{h_j-1}E_j^i, S^{h_j}E_j^i, S^{h_j+1}E_j^i, \dots, S^{h_j+s_j(i)-1}E_j^i$$
(все эти множества  не пересекаются).
Обозначив $E_{j+1}= E^1_j$, для   $i<r_j$ положим 
$$S^{h_j+s_j(i)}E_j^i = E_j^{i+1}.$$
 Набор надстроеных колонн с этого момента  рассматривается как   башня  этапа $j+1$,  состоящая из полуинтервалов  
$$E_{j+1}, SE_{j+1}, S^2 E_{j+1},\dots, S^{h_{j+1}-1}E_{j+1},$$
где  
 $$ h_{j+1}=h_jr_j +\sum_{i=1}^{r_j}s_j(i).$$

Частичное определение преобразования $S$ на этапе $j$ сохраняется на всех следующих этапах. В итоге на пространстве  $X=\cup_j X_j$ определено  обратимое преобразование $S:X\to X$, сохраняющее  стандартную меру Лебега на $X$.

Пусть конструкция $S$ ранга один обладает следующим свойством: \it пересечение 
$X_j\cap S^mX_j$ при   $h_{j}<m\leq h_{j+1}$  может содержаться 
только в одной  из колонн  $X_{i,j}$ башни $X_j$. \rm Такая конструкция 
называется  \bf сидоновской. \rm Мера пространства $X$ в этом случае бесконечна.

\vspace{2mm}
\bf  Теорема 1.  \it Пусть $T$ -- перемешивающее преобразование, $S$  -- сидоновское преобразование ранга один.  Для всякого $\eps>0$ найдется такое сопряжение $R$, что $  r(R^{-1}SR,T)<\eps$.\rm

\vspace{2mm}
\bf  Лемма 2.  \it  Пусть $T$ -- перемешивающее преобразование, $A$  -- множество конечной меры.  Фазовое пространство для всякого натурального $H$   можно представить в виде  $$X=\uu_{h>H}\uu_{i=0}^{h}T^iB_{h}, \ \ \ T^hB_{h}\cap A=\phi.$$  \rm

Доказательство леммы несложное, оставляем его  в качестве упражнения.

Доказательство теоремы кратко.  Рассмотрим часть башни Какутани:
 $$U_N=\uu_{H<h<N}\uu_{i=0}^{h}T^iB_{h},$$
которая почти поглотила множества $A_1,\dots, A_k$, фигурирующие в определении метрики $r$. Для этого выбираем большое $N$.

Возьмем очень тонкую башню $X_j$, фигурирующую в построении сидоновского преобразования $S$, намотаем   на  $U_N$ соответствующую сопряжению 
$ \tilde S=R^{-1}SR$ башню $\tilde X_j$   таким образом, чтобы выполнялись  условия:

(1) сопряженное преобразование $\tilde S$ совпадало с $T$ на $U_N$ минус  верхние  этажи $T^hB_{h}$ и минус множество сколь угодно малой меры;

(2)  после того как движущийся этаж намотанной башни покидает в $k$-тый раз верх $U_N$, он должен  путешествовать вне   $U_N$  время $s_k> 10 s_{k-1}$
($s_1> 10N$). Этот способ можно назавать сидоновской намоткой, которая 
 обеспечит нам малость меры пересечения: 
$$\mu(\tilde S^m\tilde U_N\cap U_N)\ <\ N\mu(B_j), \ \ m<h_j.$$

При $m\geq h_j$  выполняется неравенство 
$$\mu(\tilde S^m\tilde U_N\cap \tilde U_N)<\mu(\tilde U_N)/r_j$$  в силу 
определения сидоновской конструции.

Если бы выполнялось $A=A_1\cup\dots \cup A_k\subset \tilde U_N$,
то  
$$\sum_i^k 2^{-i}\left(\mu(SA_i\Delta TA_i)+\mu(S^{-1}A_i\Delta T^{-1}A_i)\right)=0$$
и  
$$\sum_{i,j=1}^k 2^{-i-j}\left|\mu(SA_i\cap A_j)-\mu(TA_i\cap A_j)\right|=0.$$
В действительности при  выборе большого $N$ и подходящей обмотке вместо равенств получим неравенства с очень малой правой частью.

Итак для  $\eps>0$ всегда найдутся 
$j$, $N$ и $R$ такие, что на пересечении  $\tilde X_j\cap U_N$ преобразования $\tilde S$ и $T$ почти совпадают, значит близки в метрике Халмоша,
но их степени  равномерно близки относительно метрики $d_w$.
Тем самым мы получим $r(R^{-1}SR,T)<\eps$ для подходящего $R$.

Множество преобразований ранга один является $G_\delta$ множеством 
(используем без особых изменений доказательство из \cite{Ba}).
Мы показали его  плотность.  Таким образом, типичность ранга один доказана.

\vspace{3mm}
Замечание.  В \cite{20} мы использовали подобный прием, но обмотки были случайными орнстейновскими и надстройки для потоков и автоморфизов ранга один также были случайными орнстейновскими.
В случае бесконечной меры возникает приятная возможность использовать  сидоновкие надстройки, что дает  простой и  конструктивный метод.   Другие результаты о сидоновских преобразованиях см. в \cite{24}.

\end{document}